\documentclass[10pt]{amsart}
\usepackage{amsmath}
\usepackage{enumerate}
\usepackage{pst-node}
\usepackage{tikz-cd}
\usepackage{amsopn}
\usepackage{amsmath, amsthm, amssymb, amscd, amsfonts, xcolor}
\usepackage[latin1]{inputenc}
\usepackage{graphicx}

\usepackage{euscript}

\usepackage{xcolor}
\usepackage{enumitem}
\usepackage{multirow}

\input xy
\xyoption{all}

\usepackage[OT2,OT1]{fontenc}
\newcommand\cyr{%
\renewcommand\rmdefault{wncyr}%
\renewcommand\sfdefault{wncyss}%
\renewcommand\encodingdefault{OT2}%
\normalfont
\selectfont}
\DeclareTextFontCommand{\textcyr}{\cyr}

\theoremstyle{plain}

\newtheorem{teo}{Theorem}[section]

\newtheorem{prop}[teo]{Proposition}

\newtheorem{rem}[teo]{Remark}

\theoremstyle{definition}

\numberwithin{equation}{section}

\title{Radial parallel sections of vector bundles}

\author{Antonio J. Di Scala}

\address{Dipartimento di Scienze Matematiche \it G.L. Lagrange \rm , Politecnico di Torino, Corso Duca degli Abruzzi 24, 10129, Torino, Italy.}

\email{antonio.discala@polito.it}
\thanks{This work was partially supported by the Simons - Foundation grant 346300 and the Polish Government MNiSW 2015-2019 matching fund.
The idea of the proof presented in this paper comes to the author whilst he was at the Banach Center at the IMPAN in Warsaw for the Simons Semester \it Symmetry and Geometric Structures. \rm The author would like to thank Mike Eastwood and Pawel Nurowski for very interesting discussions during his stay at the Simons Semester.
Remark 1.2. contains a recent discussing with Mike Eastwood.
The author is member of G.N.S.A.G.A. of I.N.d.A.M.}


\begin{document}

\begin{abstract}
We give a proof of the existence of radial (smooth) parallel sections of vector bundles endowed with a linear connection.
\end{abstract}

\maketitle

\section{Introduction}

In several talks I gave about my paper \cite{DM16} (e.g. at the meeting \it  Variet\'a reali e complesse: geometria, topologia e analisi armonica \rm, SNS, Pisa (Italy) 2014, at the Simons Semester \it Symmetry and Geometric Structures \rm at the IMPAN, Warsaw (Poland) 2017) I was asked about the proof of the following proposition:

\begin{prop}\label{main} Let $\pi:E \to M$ be a vector bundle with a connection $\nabla$. Let $(x_1,\cdots,x_n)$ be local coordinates centered at a point $p \in M$ and let $\xi_p \in E_p := \pi^{-1}(p)$. Then there exists a neighborhood $U$ of $p$ and a smooth section $\xi : U \to E$ such that $\xi(p) = \xi_p$ and \begin{equation} \label{rp} \nabla_{\frac{\partial}{\partial \rho}} \xi \equiv 0 \end{equation}
where $\frac{\partial}{\partial \rho} := \sum_{j=1}^n x_j \frac{\partial}{\partial x_j}$ is the radial vector field. That is to say, there exists a radial $\nabla$-parallel smooth section with a given initial condition $\xi_p \in E_p$.
\end{prop}

I learnt by heart the above proposition from C. Olmos when I was a PhD. student.

The standard theorem of existence and uniqueness of solutions of a system of ODE implies that the function $\xi : U \to E$ satisfying (\ref{rp})  is uniquely determined by its initial condition $\xi_p$.
People who asked me about $\xi$ were concerned about the smoothness of $\xi$ at the point $p$. The goal of this short note is to show that $\xi$ is indeed smooth at $p$.

\subsection{The proof.}

Let $k = \mathrm{rank}(E)$ and let $(\sigma_1,\cdots,\sigma_k)$ be a local smooth frame of $\pi: E \to M$ defined in a neighborhood $U$ of $p$. We assume that $U$ is contained in the domain of the coordinates $x_1,\cdots,x_n$ centered at $p$.
Let $\Gamma_{ij}^s \in C^{\infty}(U)$ be the connection symbols defined by
\[ (\nabla_{\frac{\partial}{\partial x_i}} {\sigma_j} )_{|_{p}}= \sum_{s = 1}^k \Gamma_{i j}^s(p) \cdot \sigma_s(p) \, \, .\]
For each $z \in U$ fixed let $\xi(t,z)$ be the $\nabla$-parallel section defined on the radial segment $t:[0,1] \to t \cdot z $  with initial condition $\xi_p$. So $\xi(t,z) = \sum_{j=1}^k y_j(t,z) \cdot \sigma_j(t \cdot z)$ and the coefficients $y_k(t,z)$ satisfies the ODE system:
\begin{equation} \label{ode}  \sum_{i=1}^n z_i \cdot \mathrm{M}_i \cdot \mathrm{y} + \frac{\partial \mathrm{y}}{\partial t} = 0
\end{equation}
where $\mathrm{y}$ is the column $\left[
                                 \begin{array}{c}
                                   y_1(t,z) \\
                                   \vdots \\
                                   y_k(t,z) \\
                                 \end{array}
                               \right]$ and $\mathrm{M}_i$, $(i=1,\cdots,n)$ are the $k \times k$ matrices
                               $\mathrm{M}_i := (\Gamma_{i j}^s(t \cdot z))$, $(s,j = 1,\cdots,k)$.\\

Set $f(t,\mathrm{y},z) := -\sum_{i=1}^n z_i \cdot \mathrm{M}_i \cdot \mathrm{y}$. Observe that $f(t,\mathrm{y},z)$ is a $C^{\infty}$ function of $(t,\mathrm{y},z)$ and that $\mathrm{y}(t,z)$ solves the ODE system \[ \frac{\partial \mathrm{y}}{\partial t} = f(t,\mathrm{y},z) \, .\]

According to \cite[Corollary 4.1., page 101]{H64} $\mathrm{y}(t,z)$ is a $C^{\infty}$ function of $(t,z)$. Hence $\mathrm{y}(1,z) \in C^{\infty}(U)$ and the section $\xi : U \to E$ defined as \[\xi(z) := \sum_{j=1}^k y_j(1,z) \cdot \sigma_j(z) \]
is smooth on $U$. Since $\mathrm{y}(1,t \cdot z) = \mathrm{y}(t,z)$ and $\mathrm{y}(t,z)$ satisfies (\ref{ode}) we get that $\xi$ satisfies equation (\ref{rp}). This complete the proof of our proposition. $\Box$

\subsection{Geometric interpretation} Let $\mathrm{I} = [0,1]$ be unit interval, let $\mu: I \times \mathbb{R}^n \to \mathbb{R}^n $ defined by $\mu(t,x):= t \cdot x $ and let $\pi:E \to \mathbb{R}^n$ be a vector bundle with a linear connection $\nabla$. Let $\nabla^*$ be the pull-back connection of the pull-back vector bundle $\pi^*: E^* \to \mathrm{I} \times \mathbb{R}^n$ defined by the commutative diagram
\[ \begin{tikzcd}
E^* \arrow{r}{\mu^*} \arrow[swap]{d}{\pi^*} & E \arrow{d}{\pi} \\%
\mathrm{I} \times \mathbb{R}^n \arrow{r}{\mu}& \mathbb{R}^{n}
\end{tikzcd}
\]

The usual existence, uniqueness, and smooth dependence on initial conditions of solutions to ODE applies directly to get a smooth section $\xi^*$ (of $E^*$) which is
$\nabla^*$-parallel in the direction of the factor $\mathrm{I}$ with the initial conditions given on the hyperplane $H := \{0\} \times \mathbb{R}^n$.
The restriction of $\xi^*$ to the hyperplane  $\{1\} \times \mathbb{R}^n$ gives the radial parallel section $\xi$ of $E$.

\begin{rem}Let $\mathrm{B}: \tilde{\mathbb{R}^n} \to \mathbb{R}^n$ be the blow-up at the origin and let $\tilde{\nabla}$ be the pullback connection on the pullback vector bundle $\tilde{E} \to \tilde{\mathbb{R}^n}$. The initial condition $\xi_p$ lifts to a initial condition $\tilde{\xi_p}$ on the exceptional divisor $D := \mathrm{B}^{-1}(0)$. Now the usual existence, uniqueness, and smooth dependence on initial conditions of solutions to ODE implies that there is smooth section of $\tilde{\xi}$ defined on neighborhood $\tilde{U}$ of $D$ which is $\tilde{\nabla}$-parallel along the transverse foliation of $\tilde{\mathbb{R}^n}$. In this setting Proposition \ref{main} implies that the pulldown (or pushforward) section $B^{*}(\tilde{\xi})$ of $\pi: E \to U$ is smooth at the origin and $\nabla$-parallel along the radial directions. Actually, Proposition \ref{main} is equivalent to the statement that the pulldown (or pushforward) section $B^{*}(\tilde{\xi})$ of $\pi: E \to U$ is smooth at the origin. But I could not show this directly. The problem is that, at difference as what happen in Algebraic Geometry (or under analyticity hypothesis), in the smooth setting the constancy of a regular function on the exceptional divisor does not implies that regular function is the pullback of a regular function down stairs.
\end{rem}


\end{document}